\documentclass{article}
\usepackage{amsfonts,amssymb,latexsym,amsmath,amscd,euscript}
\newcounter{num}[section]
\newcommand{\Num}{\refstepcounter{num}%
	\textbf{\arabic{section}.\arabic{num}}}

\newcommand{\Theorem}{\textbf{Theorem~}}
\newcommand{\Proof}{{\sc{Proof}}}
\newcommand{\Def}{\textbf{Definition~}}

\newcommand{\Lemma}{ \textbf{Lemma~}}

\newcommand{\Remark}{\textbf{Remark~}}
\newcommand{\Prop}{\textbf{Proposition~}}

\newcommand{\Ch}{{\mathfrak S}}
\newcommand{\GL}{{\mathrm{GL}}}

\newcommand{\Irr}{{\mathrm{Irr}}}
\newcommand{\Ind}{{\mathrm{Ind}}}

\newcommand{\al}{{\alpha}}
\newcommand{\la}{{\lambda}}

\newcommand{\Fq}{\mathbb{F}_q}

\newcommand{\Ad}{{\mathrm{Ad}}}
\newcommand{\ad}{\mathrm{ad}}

\newcommand{\gx}{{\mathfrak g}}
\newcommand{\gl}{\mathfrak{gl}}
\newcommand{\nx}{{\mathfrak n}}
\newcommand{\hx}{{\mathfrak h}}
\newcommand{\bx}{\mathfrak{b}}
\newcommand{\Bx}{{\mathfrak B}}

\newcommand{\mx}{{\mathfrak m}}
\newcommand{\px}{{\mathfrak p}}
\newcommand{\ax}{{\mathfrak{a}}}
\newcommand{\gxa}{{\mathfrak g}^a}
\newcommand{\Ux}{{\mathfrak U}}

\newcommand{\nxm}{{\mathfrak n}_-}

\newcommand{\rt}{\mathrm{right}}

\newcommand{\tGa}{{\widetilde{G^a}}}
\newcommand{\tG}{\widetilde{G}}

\renewcommand{\leq}{\leqslant}

\newcommand{\Jc}{{\mathcal J}}
\newcommand{\Kc}{{\mathcal K}}

\newcommand{\Ac}{{\mathcal A}}
\newcommand{\Bc}{{\mathcal B}}
\newcommand{\Sc}{{\mathcal{S}}}

\newcommand{\Yc}{{\mathcal Y}}

\newcommand{\Cb}{\mathbb{C}}

\newcommand{\row}{{\mathrm{row}}}
\newcommand{\col}{{\mathrm{col}}}
\newcommand{\spann}{\mathrm{span}}

\newcommand{\diag}{{\mathrm{diag}}}
\newcommand{\UT}{{\mathrm{UT}}}

\newcommand{\eps}{{\varepsilon}}

\begin{document}

\title{Supercharacter theory for the Borel contraction of the group   $\GL(n,\Fq)$}
\author{A.N.\,Panov \footnote{The paper is supported by the RFBR grant 20-01-00091a }}
\date{}
\maketitle
\begin{flushleft}
	Samara  University,
	Department of Mathematics,  443011, ul. akademika Pavlova 1, Samara, Russia
\end{flushleft}
\begin{flushleft}
	UDK 512.547, MSC 20C33 \end{flushleft}
\begin{flushright} To Aleksandr Ivanovich Generalov\\ on occasion of the  seventieth anniversary \end{flushright}

\begin{abstract}
The notion of a supercharacter theory  was proposed by P. Diaconis and I.M. Isaacs in 2008.
 A supercharacter theory for a given finite group  is a pair of the system of certain complex characters  and the partition of  group into classes that   have properties similar to the system of irreducible characters and the partition into conjugacy classes.   In the present paper, we consider the group obtained by  the Borel contraction from the  general linear group over a finite field. For this group, we construct the supercharacter theory. In terms of rook placements, we classify supercharacters and superclasses, calculate values of supercharacters on superclasses.
\end{abstract}

\section{Introduction}

The main problem in representation theory is the classification problem  of irreducible rep\-re\-sen\-tations. Nevertheless, for  some groups such as the unitriangular group over a finite field, Borel  and parabolic subgroups  and others, this problem is a extremely difficult "wild"\, problem. The same is true for the classification problem of conjugacy classes.

In the paper \cite{DI}, 2008,  P. Diaconis and I.M.Isaacs introduces the notion of a supercharacter theory. By definition, a supercharacter theory of a given finite group  $G$  is a pair  $(\Ch,\Kc)$, where $\Ch=\{\chi_1,\dots,\chi_m\}$  is a system of complex characters (representations) of the group  $G$, and  $\Kc=\{K_1, \ldots, K_m\}$ is a partition of the group $G$ such that \\
S1) the characters from $\Ch$ is pairwise orthogonal;\\
S2) each character $\chi_i$ is constant on each  $K_j$;\\
S3)  $\{1\}\in \Kc$.

 The characters from $\Ch$ are referred to as \emph{supercharacters}, and  the subsets from  $\Kc$ are  \emph{superclasses}.  Observe that  the number of supercharacters equals to the number of superclasses. One of  supercharacter theory examples  is a pair of irreducible characters and the partition of  group into conjugacy classes.

    Supercharacters are not irreducible in general; their supports (the sets of irreducible constituents)  form a partition of the  set $\Irr(G)$ of all complex irreducible characters (see \cite{DI, P0}). In its turn, superclasses  decompose into conjugate classes.

One can compare supercharacter theories as follows: a supercharacter is coarser than the other supercharacter theory if its superclasses are unions of the superclasses of  other one. Therefore, any supercharacter theory is coarser than the theory of irreducible characters.   If it is impossible to classify irreducible characters and classes of conjugacy elements, the aim is to construct a supercharacter theory with exact classification of supercharacters and superclasses that is as close to the theory of irreducible characters as possible.

The classical example of a supercharacter theory is the supercharacter theory for algebra groups over a finite field from  \cite{DI} (by definition, an algebra group
is a group of the form  $G=1+\Jc$,  where  $\Jc$ is a nilpotent associative finite dimensional algebra). In the case of the unitriangular group  $\UT(n,\Fq)$, this  supercharacter theory coincides with the theory of basic characters of C.\,	Andr\'{e} \cite{A1,A2,A3}.

In the present paper, we consider the group
$G^a$ that is obtained by the group contraction of   $GL(n,\Fq)$  in the sense of definition from   \cite{IW, GOV}. The commutator $U^a$ of the group $G^a$ is an algebra group and, therefore, it has the supercharacter theory from    \cite{DI}.
The group  $G^a$ is a group of invertible elements in the reduced associative algebra  $\gxa$, and one can construct  a supercharacter theory for this group following the paper  \cite{P1}.

The main goal of this paper is to classify supercharacters and superclasses in  $U^a$ and  $G^a$, to obtain the formula for values of supercharacters on superclasses in the exact form.

 \section{Contraction of the group  $\GL(n,\Fq)$}

Consider the group  $G=\GL(n,K)$ defined over a field  $K$ of characteristic not equal to 2.
Its Lie algebra  $\gx=\gl(n,K)$ affords decomposition   $\gl(n,K)=\nx +\hx+\nxm$, where
$\nx=\nx_+$ (respectively, $\nx_-$) is a Lie  subalgebra of strongly upper triangular  (respectively, lower triangular)  matrices, $\hx$ is the subalgebra of diagonal matrices.
The Borel subalgebra  $\bx=\hx+\nx$ coincides with the Lie algebra of triangular matrices.
The Killing form identifies  $\nx_-$ with the  dual space  $\nx^*$ of $\nx$.

Consider the  Lie algebra   $\gxa$ that is the semidirect sum of the Borel subalgebra  $\bx=\hx+\nx$ and the commutative Lie algebra  $\nx^*$,  which coincides with $\nxm$ as a linear space; the action of    $x\in\bx$ on  $\la\in\nx^*$ is the coadjoint action  $\ad^*_x(\la)$.
We refer to the Lie algebra   $\gxa$  as  \textit{the Borel contraction of the Lie algebra } $\gx$.

The group $G$ contains the triangular group $B=HN$, where $H$ is the subgroup of diagonal matrices, ~ $N$ is the unitriangular group.
Denote by   $G^a$  the  semidirect product of the Borel subgroup   $B=HN$ and the abelian  group  $N^a_-=1+\nx^*$ with multiplication  $(1+\la_1)(1+\la_2)=1+\la_1+\la_2$; the action of  $b\in B$ on $N_-^a$ is defined as   $b(1+\la)b^{-1}=1+\Ad^*_b(\la)$.
We refer to the group   $G^a$  as  \textit{the Borel contraction of the group }  $G$.

Let us show that the Lie algebra  $\gx^a$  affords the  structure of associative algebra such that $[X,Y]=XY-YX$. Indeed, in our case,  $\bx$ is an associative algebra of upper triangular matrices, $\nx^a$ is an algebra with zero multiplication.
There are  defined the left and right actions of the algebra $\bx$  (respectively, group  $B$) on   $\nx^*$ by the formulas
\begin{equation}\label{rlaction}
t\la(x)=\la(xt)\quad\mbox{and}\quad \la t(x)=\la(tx),
\end{equation}
where $t\in \bx$, ~$x\in \nx$ and $\la\in \nx^*$.
One can  verify that   $\gx^a$ is an associative algebra with respect to the following  multiplication
$$(t_1+\la_1)(t_2+\la_2)=t_1t_2+t_1\la_2+\la_1 t_2,$$
where  $t_1, t_2\in \bx$ and $\la_1,\la_2\in \nx^*$.
The radical  $\Ux =\nx+\nx^*$ of the associative algebra  $\gxa$ is a nilpotent associative subalgebra.

Observe that the group $G^a$ is the  group of  invertible elements of the associative algebra  $\gx^a$.
The group  $G^a$ contains the normal subgroup  $U^a =1+\Ux$, which is an algebra group.

Define the symmetric bilinear form on $\Ux$ as follows
\begin{equation}\label{form}
(X_1,X_2)=\la_2(x_1)+\la_1(x_2),
\end{equation}
where $X_1=x_1+\la_1$ and $X_2=x_2+\la_2$.

Easy to show that
\begin{equation}\label{uniform}
(gX_1,X_2) = (X_1,X_2g)
\end{equation}
for all  $g\in G^a$ and  $X_1,X_2\in \Ux$.
The bilinear form  (\ref{form}) is nondegenerate;  this enables us to identify  $\Ux$ with its dual space $\Ux^*$.  In doing so, the adjoint orbits in $\Ux$ are identified with the coadjoint ones in $\Ux^*$

We  call a \textit{root} any integer pair  $(i,j)$, where $i\ne j$ and $1\leq i, j\leq n$.
This definition is common in the  theory of basic characters and suitable in  dealing with rook spacements, right stabilizers and others. The pair  $\alpha=(i,j)$ is assigned to a root $\alpha=\epsilon_i - \epsilon_j$ in the Lie theory.

For  each root $\al=(i,j)$, we call  $i$  \emph{the row number} of   $\al$, and  $j$  \emph{ the column number} of $\al$ (notations $i=\row(\al)$ and $j=\col(\al)$). On the set of root  $\Delta$ one can define the partial addition  $(i,j)+(j,k)=(i,k)$. For each root $\al=(i,j)$, we define the \emph{opposite root} $\al^t=(j,i)$.

We say that  $\al=(i,j)$ is  \textit{a positive root} (respectively, negative root)  if  $i<j$ (respectively,   $i>j$). The set of all roots  $\Delta$ is a union of the subset of positive roots  $\Delta_+$ and the subset of negative roots $\Delta_-$. Attach to each positive root  $\al=(i,j)$,~ $i<j$, the matrix unit  $E_\al=E_{ij}\in \nx$. Respectively, attach to each negative root  $\beta=(i,j)$,~ $i>j$, the matrix unit $F_\beta=F_{ij}\in \nx_-=\nx^*$. Observe that  $(E_{ij},F_{km})=\delta(j,k)\delta(i,m)$.

We realize the associative algebra  $\Ux=\nx+\nx^*$ as a subspace spanned by the union of two systems of matrix units  $\{E_{ij}:~ i<j\}$ and  $\{ F_{ij}: i>j\}$.
The structure relations have the form
$$ E_{ij}E_{km}=\delta(j,k) E_{im},\quad\quad F_{ij}F_{km}=0,$$
$$ E_{ij}F_{km}=\left\{\begin{array}{cl} F_{im}&,~\mbox{if}~~j=k~~\mbox{and}~~i>m,\\
0&,~\mbox{in~other~ cases}.\end{array}\right.
$$
$$ F_{ij}E_{km}=\left\{\begin{array}{cl} F_{im}&,~\mbox{if}~~j=k~~\mbox{and}~~i>m,\\
0&,~\mbox{in~other~ cases}.\end{array}\right.
$$
\Def\Num. We refer to a  subset $D\subset \Delta$ as a rook placement  if there is  at most one root of $D$ in any row and  column.
The subset  $D$ decomposes $D=D_+\cup D_-$, where $D_+=D\cap \Delta_+$ and $D_-=D\cap\Delta_-$.

Observe that here the number of rooks may be smaller than the size of board $n$.
Let  $\phi$  be a map  $D\to\Fq^*$. The map  $\phi$  is uniquely defined by its restrictions  $\phi_+$ and $\phi_-$ on $D_+$ and $D_-$ respectively.
For any pair $(D,\phi)$, we define the element  $ u_{D,\phi}=1+X_{D,\phi}\in U^a$, where
\begin{equation}\label{elemDphi}
X_{D,\phi}=\sum_{\gamma_+\in D_+} \phi(\gamma_+) E_{\gamma_+} + \sum_{\gamma_-\in D_-} \phi(\gamma_-) F_{\gamma_-}.
\end{equation}

\section{Supercharacter theory for the group  $U^a$}

In what follows,  $K$ is a finite field   $\Fq$ of characteristic not equal to  2.
Since the group   $U^a$  is an algebra group, the supercharacter theory of  P. Diaconis and  I.M. Isaacs is defined for it \cite{DI}.

\subsection{Supercharacter theory for algebra groups}\label{algr}
Let   $G=1+\Jc$  be an algbera group, where  $\Jc$ is an associative finite dimensional  nilpotent algbera over  the finite field $\Fq$.   The group $G$ acts on  $J$ by the left and right multiplication. Superclasses of  $G$ are the subsets of the form  $K(x)=1+GxG$.

By the formulas  (\ref{rlaction}), we define the left and right actions of   $G$ on the dual space  $\Jc^*$. Respectively, we define the action of the group  $G\times G$ on  $\Jc^*$ by the left-right multiplication.  The dual space $\Jc^*$ is divided into  $G\times G$ orbits. For any  $\la\in\Jc^*$, consider the stabilizer  $G_{\la,\rt}$ of the right action of  $G$ on $\Jc^*$.  The subgroup  $G_{\la,\rt}$ is an algebra group
$G_{\la,\rt} =1+ \Jc_{\la,\rt}$, where
$$\Jc_{\la,\rt}=\{x\in\Jc:~ \la (xy)=0~~ \mbox{for ~ every}~ y\in\Jc\}.$$
Fix  a nontrivial character  $t\to\eps^t$ of the additive group of the field  $\Fq$ with values in   $\Cb^*$. The  formula  $\xi_\la(1+x) =\eps^{\la(x)}$ defines the linear character (one dimensional representation) of  subgroup   $G_{\la,\rt}$. Consider the  character $\chi_\la$ of $G$ induced from the character  $\xi_\la$ of $G_{\la,\rt}$ to $G$.
\\
\Theorem \Num \label{thalgr} \cite{DI}. \emph{The  systems of subsets  $\{K(x)\}$ and characters  $\{\chi_\la\}$, where $x$ and  $\la$ run through the systems of representatives of  $(G\times G)$-orbits in $J$ and  $J^*$ respectively, give rise to a supercharacter theory for the algebra group  $G$.}\\
\Remark~\Num\label{remalgr}. In definition of the character  $\chi_\la$, one can replace the stabilizer of  right action by the stabilizer  of left action (see \cite{DI}).

For the group  $N=\UT(n,\Fq)$, this supercharacter theory coincides with the theory of basic characters of C. Andr\'{e} and afford complete description in terms of rook placements. The algebra
$\Jc$  equals to  $\nx$ and the dual space  $\Jc^*$ equals to  $\nx^*=\nx_-$.
\\
\Prop\Num\label{basic} \cite{DI,A2,A3}.
\emph{1. In each   $(N\times N)$-orbit in  $\nx$,  there is a unique element of the form  $X_{D_+,\phi_+}$.
2. In each  $(N\times N)$-orbit in  $\nx^*$,  there is a unique element of the form  $X_{D_-,\phi_-}$.}

\subsection{Superclasses in  $U^a$}

The group  $U^a=1+\Ux$ is an algebra group.
The superclass of  element $u=1+X$ from  $U^a$ is defined as  $K(u)=1+U^a XU^a$ (see previous subsection). The next theorem provides a classification of superclasses.\\
\Theorem\Num\label{supclass}.\emph{ 1. Each superclass containes an element  $ u_{D,\phi}$  for some rook placement  $D$ in the root system $\Delta$ and a map $\phi:\Delta\to\Fq$.\\
2. Two elements $ u_{D,\phi}$ and  $ u_{D',\phi'}$ belong  to a common class if and only if  $(D,\phi)=(D',\phi')$.}\
 \Proof. \textit{Item 1}.
Let  $X=x+\la\in \Ux$,  where $x\in \nx$ and  $\la\in\nx^*$. Let us show that there exists  $X_{D,\phi}$ belonging to  $U^a XU^a$. The proof is divided into several sub-items.\\
1) There exists   $a,b\in N$ such that  $axb=X_{D_+,\phi_+}$ for some rook placement  $D_+$ in $\Delta_+$ and the map  $\phi_+:D_+\to \Fq^*$ (see Proposition 3.3)).  	
We take  $x_+=X_{D_+,\phi_+}$. Then $aXb=x_++a\la b$.\\
2) It is sufficient to prove the statement of Item 1 for element of the form
$X=x_++\la$. In this sub-item, we consider decompositions of  $\nx$ and $\nx^*$ associated with $D_+$ and show that  the $(U^a\times U^a)$-orbit of  $X$ contains some element of the form $x_++\la$  obeying the condition   $\la\in \ell_0^\perp\cap r_0^\perp$ (see formula  (\ref{nstar})).

Consider the decomposition  $\nx=\ell_+\oplus \ell_0$ of $\nx$ into the sum of left ideals
$$\ell_+= \spann \{E_{ij}:~ j\notin \row(D_+)\}=\{y\in\nx:~ yx_+=0\},\quad\quad
\ell_0=\spann \{E_{ij}: ~ j\in \row(D_+)\}.$$
We obtain $\nx^*=\ell_+^\perp\oplus\ell_0^\perp$, where
$$\ell_+^\perp=\spann\{F_{km}:~ k\in \row(D_+)\}=x_+\nx^*,\quad\quad
\ell_0^\perp=\spann\{F_{km}:~ k\notin \row(D_+)\}.$$
We have $\nx^*=x_+\nx^*\oplus\ell_0^\perp$.\\
{\bf Example}. $n=4$, ~$D=\{\gamma=(2,3)\}$. Then
{\small
	$$x_+= \left(\begin{array}{cccc}
	0&0&0&0\\
	0&0&1&0\\
	0&0&0&0\\
	0&0&0&0\\
	\end{array}\right),\qquad \ell_+=\left\{\left(\begin{array}{cccc}
	0&0&*&*\\
	0&0&*&*\\
	0&0&0&*\\
	0&0&0&0\\
	\end{array}\right)\right\},\qquad  \ell_0=\left\{\left(\begin{array}{cccc}
	0&*&0&0\\
	0&0&0&0\\
	0&0&0&0\\
	0&0&0&0\\
	\end{array}\right)\right\}
	$$
	$$ \ell^\perp_+ = x_+\nx^* = \left\{\left(\begin{array}{cccc}
	0&0&0&0\\
	*&0&0&0\\
	0&0&0&0\\
	0&0&0&0\\
	\end{array}\right)\right\},\qquad  \ell^\perp_0=\left\{\left(\begin{array}{cccc}
	0&0&0&0\\
	0&0&0&0\\
	*&*&0&0\\
	*&*&*&0\\
	\end{array}\right)\right\}.
	$$
}
Analogously, we consider the decomposition  $\nx=r_+\oplus r_0$ of $\nx$ into the sum of  right ideals
$$r_+= \spann \{E_{ij}:~ i\notin \col  (D_+)\}=\{y\in\nx:~ x_+y=0\},\quad\quad
r_0=\spann \{E_{ij}: ~ i\in \col (D_+)\}.$$
We obtain $\nx^*=r_+^\perp\oplus r_0^\perp$, where
$$r_+^\perp=\spann\{F_{km}:~ m\in \col(D_+)\}=\nx^*x_+,\quad\quad
r_0^\perp=\spann\{F_{km}:~ m\notin \col(D_+)\}.$$
We have $\nx^*=\nx^*x_+\oplus r_0^\perp$.

{\bf Example}. $n=4$, ~$D=\{\gamma=(2,3)\}$. Then
{\small
	$$x_+= \left(\begin{array}{cccc}
	0&0&0&0\\
	0&0&1&0\\
	0&0&0&0\\
	0&0&0&0\\
	\end{array}\right),\qquad r_+=\left\{\left(\begin{array}{cccc}
	0&*&*&*\\
	0&0&*&*\\
	0&0&0&0\\
	0&0&0&0\\
	\end{array}\right)\right\},\qquad  r_0=\left\{\left(\begin{array}{cccc}
	0&0&0&0\\
	0&0&0&0\\
	0&0&0&*\\
	0&0&0&0\\
	\end{array}\right)\right\}
	$$
	$$ r^\perp_+ = \nx^*x_+ = \left\{\left(\begin{array}{cccc}
	0&0&0&0\\
	0&0&0&0\\
	0&0&0&0\\
	0&0&*&0\\
	\end{array}\right)\right\},\qquad  r^\perp_0=\left\{\left(\begin{array}{cccc}
	0&0&0&0\\
	*&0&0&0\\
	*&*&0&0\\
	*&*&0&0\\
	\end{array}\right)\right\}.
	$$
}

Since the subspaces  $\ell_+$, $\ell_0$, $r_+$ and $r_0$ are spanned by the elements of the standard basis $\{E_{ij},F_{km}\}$,  we get
\begin{equation}\label{nstar}
\nx^*=\left(x_+\nx^* + \nx^*x_+\right)\oplus \left(\ell_0^\perp\cap r_0^\perp\right).
\end{equation}

{\bf Example}. $n=4$, ~$D=\{\gamma=(2,3)\}$. Then
{\small
	$$x_+= \left(\begin{array}{cccc}
	0&0&0&0\\
	0&0&1&0\\
	0&0&0&0\\
	0&0&0&0\\
	\end{array}\right),\quad\quad  x_+\nx^* + \nx^*x_+ =\left\{\left(\begin{array}{cccc}
	0&0&0&0\\
	*&0&0&0\\
	0&0&0&0\\
	0&0&*&0\\
	\end{array}\right)\right\},$$
	$$\ell^\perp_0\cap r^\perp_0=\left\{\left(\begin{array}{cccc}
	0&0&0&0\\
	0&0&0&0\\
	*&*&0&0\\
	*&*&0&0\\
	\end{array}\right)\right\}.
	$$
}

Let    $X=x_++\la\in \Ux$. The element $\la$ decomposes   $\la=x_+\nu_1+\nu_2x_++\la_0$, where $\la_0\in \ell_0^\perp\cap r_0^\perp$, ~ $\nu_1,\nu_2\in \nx^*$.
We have
$$(1-\nu_2)X(1-\nu_1) = (1-\nu_2)(x_++\la)(1-\nu_1)= x_+ - x_+\nu_1 -\nu_2x_++\la= x_++\la_0.$$
3) So, we assume that   $X=x_++\la_0$, where $\la_0(\ell_0)=\la_0(r_0)=0$.  Let us show that the  $(U^a\times U^a)$-orbit of   $X$ contains an element of the  form   $x_+ +\la'$, where $$\la'=X_{D_-',\phi_-'}$$ for some rook placement  $D_-'$ in $\Delta_-$ and $\phi_-':D_-'\to\Fq^*$.

For $\la_0$, there exist  $a,b\in N$ such that  $a\la_0b=\la'$, where $\la'$ have the mentioned form    (see Proposition \ref{basic}). It remains to prove that one can choose  $a,b$ such that  $ax_+b=x_+$.

Consider the algebra subgroups  $L_+=1+\ell_+=\{a_+\in N:~ a_+x_+=x_+\}$ and $L_0=1+\ell_0$.
Show that  $a_0\la_0=\la_0$ for any $a_0\in L_0$. Indeed, let  $a_0=1+y_0$,~ $y_0\in\ell_0$. Since  $\ell_0$ is a left ideal  and  $\la_0(\ell_0)=0$, for any  $x\in\nx$, we have $a_0(\la_0)(x)=\la_0(xa_0)=\la_0(x+xy_0)=\la_0(x)$.

Observe that  any element  $a\in N$  is uniquely presented in the form  $a=a_+a_0$, ~ $a_+\in L_+$ and $a_0\in L_0$. We obtain
$a\la_0=a_+a_0(\la_0)=a_+\la_0$, where $a_+x_+=x_+$.

Similarly, we define the subgroups  $R_+=1+r_+=\{b_+\in N:~ x_+b_+=x_+\}$ and $R_0=1+r_0$.  Then  $\la_0 b_0=\la_0$ for any $b_0\in R_0$.
Any element  $b\in N$ is uniquely presented in the form
$b=b_0b_+$,~ $b_+\in R_+$ and $b_0\in R_0$).
We have $\la_0b=\la_0b_0b_+=\la_0b_+$, where $x_+b_+=x_+$.

Finally, we get  $\la'=a\la_0b=a_+\la_0b_+$, where $a_+x_+b_+=x_+$. Then $a_+(x_++\la)b_+=x_++\la'$. \\
4) It is proved above that the  $(U^a\times U^a)$-orbit of $X$ contains the element  $X'=x_++\la'$, where $x_+=X_{D_+,\phi_+}$ and $\la'=X_{D_-',\phi_-'}$.
Finally, as in  3),  we  replace  $X'$ by a suitable element  of the form   $(1+\nu_1)X'(1+\nu_2)$ deleting  rooks from   $D_-'$ belonging to rows from  $\row(D_+)$  and columns from $\col(D_+)$.
We get the element   $X_{D,\phi}$, where $D=(D_+,D_-)$ is a rook placement from  $\Delta$ and  $\phi$ is a map  $D\to\Fq^*$, which coincides with  $\phi_+$ on  $D_+$  and with $\phi_-'$ on $D_-$. This proves the statement of Item 1.\\
\textit{Item 2.} Suppose $ u_{D,\phi}$ and $ u_{D',\phi'}$ belong to a common superclass. We aim to show these elements  coincide. There exist  $A,B\in U^a$ such that  $X_{D',\phi'}=AX_{D,\phi}B$. Present the elements $A$ and  $B$ as products
$A=(1+\nu_1)a$ and $B=b(1+\nu_2)$, where $a,b\in N$ and $\nu_1,\nu_2\in \nx^*$.
Then $$X_{D'_+,\phi'_+} =X_{D',\phi'}\bmod N_-^a=AX_{D,\phi}B\bmod N_-^a=
aX_{D_+,\phi_+}b.$$
Proposition  \ref{basic} implies  $(D'_+,\phi'_+)=(D'_+,\phi_+)$. Denote $x_+=X_{D_+,\phi_+}$. Then $ax_+b=x_+$.

We obtain  $X_{D,\phi}=x_++\la$ and $X_{D',\phi'}=x_++\la'$, where $\la=X_{D_-,\phi_-}$ and $\la'=X_{D'_-,\phi'_-}$.
We have $x_++\la'= (1+\nu_1)a(x_++\la)b(1+\nu_2)=x_++x_+\nu_2+\nu_1x_+ +
a\la b$ ~ and

\begin{equation}\label{xplus}
\la'= x_+\nu_2+\nu_1x_+ +
a\la b.
\end{equation}

Consider the decomposition  $\nx=\mx_+\oplus\mx_0$ of $\nx$ into the sum of right ideals  $$\mx_+=\spann\{E_{ij}:~ i\notin\row(D_+)\}=\{x\in\nx:~\ell_0x=0\},\quad\quad
\mx_0=\spann\{E_{ij}:~ i\in\row(D_+)\}.$$
{\bf Example}. $n=4$, ~$D=\{\gamma=(2,3)\}$. Then
{\small
	$$x_+= \left(\begin{array}{cccc}
	0&0&0&0\\
	0&0&1&0\\
	0&0&0&0\\
	0&0&0&0\\
	\end{array}\right),\qquad \mx_+=\left\{\left(\begin{array}{cccc}
	0&*&*&*\\
	0&0&0&0\\
	0&0&0&*\\
	0&0&0&0\\
	\end{array}\right)\right\},\qquad  \mx_0=\left\{\left(\begin{array}{cccc}
	0&0&0&0\\
	0&0&*&*\\
	0&0&0&0\\
	0&0&0&0\\
	\end{array}\right)\right\}
	$$}
Define the algebra subgroups  $M_+=1+\mx_+$ and $M_0=1+\mx_0$.
Each element  $a\in N$ is presented in the form  $a=c_0c_+$, where $c_0\in M_0$ and
$c_+\in M_+$.
Observe that
$$\mx_0\nx^*=\spann\{F_{ij}:~i\in\row(D_+)\} =  x_+\nx^*.$$ We obtain $a\la=c_0c_+\la=c_+\la\bmod x_+\nx^*$.

Similarly,  $\nx$ decomposes into the  sum of left ideals $\nx=\px_+\oplus\px_0$, where $$\px_+=\spann\{E_{ij}:~ j\notin\col(D_+)\}=\{x\in\nx:~xr_0=0\},\quad\quad
\px_0=\spann\{E_{ij}:~ j\in\col(D_+)\}.$$
{\bf Example}. $n=4$, ~$D=\{\gamma=(2,3)\}$. Then
{\small
	$$x_+= \left(\begin{array}{cccc}
	0&0&0&0\\
	0&0&1&0\\
	0&0&0&0\\
	0&0&0&0\\
	\end{array}\right),\qquad \px_+=\left\{\left(\begin{array}{cccc}
	0&*&0&*\\
	0&0&0&*\\
	0&0&0&*\\
	0&0&0&0\\
	\end{array}\right)\right\},\qquad  \px_0=\left\{\left(\begin{array}{cccc}
	0&0&*&0\\
	0&0&*&0\\
	0&0&0&0\\
	0&0&0&0\\
	\end{array}\right)\right\}
	$$}
There are defined two algebra subgroups  $P_+=1+\mx_+$ and $P_0=1+\mx_0$.
Each element $b\in N$ is  presented in the form  $b=d_+d_0$, where $d_0\in P_0$ and
$d_+\in P_+$.
Observe that
$$\nx^*\px_0=\spann\{F_{ij}:~j\in\col(D_+)\} =  \nx^*x_+.$$
Then $\la b=\la d_+d_0=\la d_+\bmod \nx^*x_+$.

It follows from  (\ref{xplus}) that
\begin{equation}\label{clad}
\la' = c_+\la d_+\bmod(x_+\nx^*+\nx^*x_+).
\end{equation}
By definition  $\la,\la'\in \ell_0^\perp\cap r_0^\perp$. Let us show that
$ c_+\la d_+\in \ell_0^\perp\cap r_0^\perp$.
Indeed, for any  $y_0\in\ell_0$, we have  $c_+\la d_+(y_0) =
\la(d_+y_0c_+)$. Since  $\ell_0$ is a left ideal and $\ell_0\mx_+=0$, we get  $d_+y_0c_+=d_+y_0\in\ell_0$.
Therefore $c_+\la d_+(y_0)=\la(d_+y_0) =0$.
Similarly,   $c_+\la d_+(z_0)=0$ for any  $z_0\in r_0$. This proves
$c_+\la d_+\in \ell_0^\perp\cap r_0^\perp$.

The elements  $\la'$ and  $c_+\la d_+$ belong to  $ \ell_0^\perp\cap r_0^\perp$. It follows from   (\ref{nstar}) and  (\ref{clad}) that  $\la' = c_+\la d_+$. Proposition 3.3 implies   $(D'_-,\phi'_-) = (D_-,\phi_-) $. $\Box$

\subsection{Supercharacters of the group  $U^a$}\label{subsecsupU}

We begin from the classification of  $(U^a\times U^a)$-orbits in  $\Ux^*$ (see subsection  \ref{algr}).
For the pair  $(D,\phi)$, where $D$ is a rook placement in $\Delta$ and $\phi: D\to\Fq^*$, we consider  $\Lambda_{D,\phi}\in\Ux^*$ such that  $$\Lambda_{D,\phi}(Y)=(X_{D^t,\phi^t}, Y),$$
where $\phi^t(\gamma) =\phi(\gamma^t)$ for any $\gamma\in D$.  Theorem  \ref{supclass} implies a classification of  $(U^a\times U^a)$-orbits in $\Ux^*$.\\
\Theorem\Num\label{supUU}.\emph{ 1.   Each $(U^a\times U^a)$-orbit in  $\Ux^*$ contains an element  $\Lambda_{D,\phi}$  for some  rook placement  $D$ in the root system  $\Delta$ and a map $\phi:D\to\Fq^*$.\\
2. Two elements  $\Lambda_{D,\phi}$ and $\Lambda_{D',\phi'}$ belong to a  common  $(U^a\times U^a)$-orbit if and only if  $(D,\phi)=(D',\phi')$.}

We aim to give a description of the  stabilizer of   $\Lambda_{D,\phi}$ with respect to the right (respectively, left) action of  the group  $U^a$ on $\Ux^*$ (for the case  $D=D_+$, see  \cite{A4,VERY}).   For each  $\gamma=(i,j)\in \Delta$ consider  the subset of roots  $\Sc(\gamma)=\Sc_+(\gamma)\cup \Sc_-(\gamma)$, where $\Sc_\pm(\gamma)\subset \Delta_\pm$ is defined as follows:
\begin{enumerate}
	\item
	If  $\gamma>0$, then $\Sc_-(\gamma)=\varnothing$ and $\Sc_+(\gamma)$ consists of the roots  $\al=(i,k)$, ~ $i<k<j$.
\item
If $\gamma<0$, then $\Sc_-(\gamma)$ consists of the roots  $\al=(i,k)$, ~ $1\leq k <j$.
\item
 If $\gamma<0$, then $\Sc_+(\gamma)$ consists of the roots $\al=(i,k)$, ~ $i< k $.

\end{enumerate}

If  $\al=(i,k)\in S(\gamma)$, then we denote  $\al^\star=(k,j)$.
Respectively,   $\Sc^\star(\gamma)$ is the set of roots  $\al^\star$, where   $\al\in\Sc(\gamma)$.  Observe that in the above-enumerated cases  1), 2) and 3) the following equalities hold  $E_\al E_{\al^\star}=E_\gamma$, ~ $F_\al E_{\al^\star}=F_\gamma$ and $E_\al F_{\al^\star}=F_\gamma$ respectively.\\
{\bf Example}.
We depict the roots $\gamma\in D$ on $n\times n$-table by the symbols  $\otimes$, roots from  $\Sc(\gamma)$ and $\Sc^\star(\gamma)$ by symbols  "$+$"\, and  " $-$"\, respectively. We fill  diagonal squares with numbers of corresponding  rows (columns). The other squares we fill with  dots.
For example, for $n=7$ and roots $\gamma_1=(2,5)$ and $\gamma_2=(5,3)$ the tables have the form

{\small
	$$
	\left(\begin{array}{ccccccc}
	1 &\cdot&\cdot&\cdot&\cdot&\cdot&\cdot\\
	\cdot &2&+&+&\otimes&\cdot&\cdot\\
	\cdot &\cdot&3&\cdot&-&\cdot&\cdot\\
	\cdot &\cdot&\cdot&4&-&\cdot&\cdot\\
	\cdot &\cdot&\cdot&\cdot&5&\cdot&\cdot\\
	\cdot &\cdot&\cdot&\cdot&\cdot&6&\cdot\\
	\cdot&\cdot&\cdot&\cdot&\cdot&\cdot&7\\
	\end{array}\right)\hspace{3cm}
	\left(\begin{array}{ccccccc}
	1 &\cdot&-&\cdot&\cdot&\cdot&\cdot\\
	\cdot &2&-&\cdot&\cdot&\cdot&\cdot\\
	\cdot &\cdot&3&\cdot&\cdot&\cdot&\cdot\\
	\cdot &\cdot&\cdot&4&\cdot&\cdot&\cdot\\
	+ &+&\otimes&\cdot&5&+&+\\
	\cdot &\cdot&-&\cdot&\cdot&6&\cdot\\
	\cdot&\cdot&-&\cdot&\cdot&\cdot&7\\
	\end{array}\right)$$}
Here $\Sc_+(\gamma_1)=\{(2,3), (2,4)\}$, ~ $\Sc_-(\gamma_1)=\varnothing$, ~ $\Sc_+(\gamma_2)=\{(5,6), (5,7)\}$, ~ $\Sc_-(\gamma_2)=\{(5,1), (5,2)\}$.

Let  $\Ux_\gamma$  be a subspace spanned by the matrix units  $E_\beta$, where $\beta\notin \Sc_+(\gamma)$, and  $F_\beta$, where $\beta\notin \Sc_-(\gamma)$. The subspace   $\Ux_\gamma$ is a subalgebra (more precisely, right ideal) in $\Ux$; respectively,   $U_\gamma=1+\Ux_\gamma$ is an algebra subgroup in  $U^a$.
Analogously,  $\Sc^\star(D)$ defines the subalgebra (left ideal) in $\Ux^\star_\gamma$ and the algebra subgroup  $U^\star_\gamma=1+\Ux^\star_\gamma$.

Take $$\Sc(D)=\bigcup_{\gamma\in D} \Sc(\gamma), \hspace{1cm} \Sc_\pm(D)=\bigcup_{\gamma\in D} \Sc_\pm(\gamma), \hspace{1cm} \Sc^\star(D)=\bigcup_{\gamma\in D} \Sc^\star(\gamma),$$
$$ \Ux_D= \bigcap_{\gamma\in D} \Ux_\gamma, \hspace{2cm} \Ux^\star_D= \bigcap_{\gamma\in D} \Ux^\star_\gamma.$$
Similarly, we define the algebra subgroups  $U_D=1+\Ux_D$ and $U^\star_D=1+\Ux^\star_D$.
By direst calculations  we verify the following statement.
\\
\Lemma\Num. \emph{The subgroup  $U_D$ (respectively, $U'_D$) is the stabilizer of $\Lambda_{D,\phi}$ with respect to the right (respectively, left) action of  $U^a$ on $\Ux^*$.}

The formula  $$\xi_{D,\phi}(u)=\eps^{\Lambda_{D,\phi}(X)}, \quad u=1+X, \quad X\in \Ux_D.$$
defines the character of the subgroup  $U_D$.
Consider the character  $$\chi_{D,\phi}=\Ind(\xi_{D,\phi}, U_D,U^a).$$

Observe that the character  $\chi_{D,\phi}$ equals to the one  induced from the  character  $\xi_{D,\phi}$ of subgroup $U^\star_D$.
The below theorem follows from Theorem \ref{thalgr}.
\\
\Theorem\Num\label{thDphi}. \emph{The systems of characters  $\{\chi_{D,\phi}\}$ and classes  $\{K_{D,\Phi}\}$ give rise to a supercharacter theory of  $U^a$.}

Our next goal is to obtain a formula for supercharacter values on superclasses. In the case $D=D_+$, the below formula coincides with corresponding formulas  from the papers  \cite{A3,VERY}.

 We need the following notations.  Let $D$ and $D'$ be two rook placements in $\Delta$. Consider the subset  	$R(D,D')$ that consists of  all elements    $\al\in \Sc(D)$ obeying  one of the following conditions:
 \begin{enumerate}
	\item If $\al\in \Sc_+(D)$, then there exists $\beta\in D'_+$  such that  $\al+\beta\in S(D)$.
	\item If $\al\in \Sc_+(D)$, then there exists $\beta\in D'_-$   such that  $\al+\beta\in S(D)$.
	\item If $\al\in \Sc_-(D)$, then there exists $\beta\in D'_+$ such that   $\al+\beta\in S(D)$.
\end{enumerate}
Denote by $r(D,D')$ the number of elements in $R(D,D')$. \\
{\bf Example}. Let $D=\{(2,6), (5,3)\}$ and $D'=\{(1,2), (3,5), (6,7), (7,1)\}$.
  Mark the roots  from  $D$, ~$D'$ and $\Sc(D)$ by symbols $\otimes$, ~ $\Box$~ and ~$+$ respectively.
{\small
	$$
	\left(\begin{array}{ccccccc}
	1 &\Box&\cdot&\cdot&\cdot&\cdot&\cdot\\
	\cdot &2&+&+&+&\otimes&\cdot\\
	\cdot &\cdot&3&\cdot&\Box&\cdot&\cdot\\
	\cdot &\cdot&\cdot&4&\cdot&\cdot&\cdot\\
	+ &+&\otimes&\cdot&5&+&+\\
	\cdot &\cdot&\cdot&\cdot&\cdot&6&\Box\\
	\Box&\cdot&\cdot&\cdot&\cdot&\cdot&7\\
	\end{array}\right)$$}
In the example,  $R(D,D')= \{(2,3), (5,1), (5,6), (5,7) \}$ and $r(D,D')=4$.

Define  $$\delta(D,D')=\left\{\begin{array}{cc} 1,&\mbox{if} ~~ D'\cap \left(\Sc(D)\cup \Sc'(D)\right)=\varnothing,\\
0,&\mbox{otherwise;}\end{array}\right.$$
$$ c(\phi,\phi')= \sum_{\gamma\in D\cap D'}\phi(\gamma)\phi'(\gamma),$$
$$ s(D)=|\Sc(D)|, \qquad m(D,D') = s(D)-r(D,D').$$

\Theorem\Num\label{Uvalue}. \emph{The value of the supercharacter  $\chi=\chi_{D,\phi}$ on the superclass
$K=K_{D',\phi'}$ equals to }
$$\chi(K)=\delta(D,D')q^{m(D,D')}\eps^{c(\phi,\phi')}.$$
\Proof.  Denote $X=X_{D',\phi'}$ and $u=1+X$.  Since supercharacters are constant on superclasses, we have $\chi(K)=\chi(u)$. Consider the subset  $S(D)$ in  $U^a$ consisting of elements of the form
\begin{equation}\label{SEF}
1+\sum_{\al\in \Sc_+(D)} a_\al E_\al + \sum_{\al\in \Sc_-(D)}a_\al F_\al,
\end{equation}
where $a_\al \in \Fq$. Easy to show that any element  $u$ from the group $U^a$ is uniquely presented in the form $u=sv$, where $s\in S(D)$ and $v\in U_D$.
The character  $\chi$ is calculated by the standard formula
\begin{equation}\label{chistandard}
\chi(u)= \sum_{s\in S(D)}\dot{\xi}(sus^{-1}),
\end{equation}
where $\dot{\xi}(u)$ equals to $\xi_{D,\phi}(u)$ on $U_D$ and zero outside  $U_D$. \\
\emph{Item 1}.
Let $\beta\in D'\cap \Sc(D)$. Let us show that  $\chi(K)=0$.
Since  $sus^{-1}= 1+sXs^{-1}$, it is sufficient to prove  $sXs^{-1} \notin \Ux_D$ for any  $s\in S(D)$. The subalgebra $\Ux_D$  is a right ideal in   $\Ux$; the element $sXs^{-1} \in \Ux_D$ if and only if  $sX\in\Ux_D$.

 Suppose $\beta\in D'\cap \Sc(D)$.
If $\beta=(i,j)\in D_+'$, then the $j$th column  $(sX)_j$ of the matrix  $sX$ has the form
$$(sX)_j = \phi'(\gamma)E_{ij}+\sum_{k<i}a_{kj}E_{kj}\bmod \nx_-^a.$$
If $\beta=(i,j)\in D_-'$, then  $$(sX)_j = \phi'(\gamma)F_{ij}+\sum_{k<i}a_{kj}F_{kj}.$$
In both cases   $sX\notin\Ux_D$, this prove the statement of Item 1.\\
\emph{Item 2.}  Suppose $ D'\cap \Sc^\star(D)\ne \varnothing$. Let us show  that  $\chi(K)=0$.
According to Remark  \ref{remalgr}, $\chi_{D,\phi}$ is the character induced form the character $\xi_{D,\phi}$ of left stabilizer  $U^\star_D$. The statement is proved analogously to the proof of Item 1. \\
\emph{Item 3.}  Suppose  $D'\cap \left(\Sc(D)\cup \Sc^\star(D)\right)=\varnothing$.
Then for any  $\beta \in D'$ the elements  $E_\beta$ and $F_\beta$ belong to   $\Ux_D$. Therefore   $X\in \Ux_D$.
For $s\in S(D)$ of the form (\ref{SEF}) and $X$ of the form (\ref{elemDphi}), we obtain
\begin{equation}\label{sX}
sX=(1+\sum_{\al\in \Sc_+(D)} a_\al E_\al + \sum_{\al\in \Sc_-(D)}a_\al F_\al)(\sum_{\beta\in D_+'} \phi'(\beta) E_{\beta} + \sum_{\beta\in D'_-} \phi'(\beta) F_\beta ).
\end{equation}

After expansion of  (\ref{sX}), we notice there is no like terms. Indeed, if  $E_{\al_1} E_{\beta_1}=E_{\al_2} E_{\beta_2}\ne 0$, then $\al_1+\beta_1=\al_2+\beta_2$ and $\col(\beta_1)=\col(\beta_2)$. As there is no two elements in $D'$ lying in the same column, we have $\beta_1=\beta_2$ and, therefore,  $\al_1=\al_2$.

If  $\al\in R(D,D')$, then there exists  the element $\beta\in D'$ such that the product of corresponding matrix units (equals to $E_\al E_\beta$, ~$E_\al F_\beta$, or $F_\al E_\beta$) does not belong to $\Ux_D$.
It follows if  $a_\al\ne 0$ for some   $\al\in R(D,D')$, then $sX\notin \Ux_D$. Since $\Ux_D$ is a right ideal, we have $sXs^{-1}\notin \Ux_D$. Then  $\dot{\xi}(sus^{-1})=0$.

 Consider the case  $a_\al=0$ for each   $\al\in R(D,D')$. 	If $\al\notin R(D,D')$, then for each  $\beta\in D' $ the elements  $E_\al E_\beta$, ~ $E_\al F_\beta$,~  $F_\al E_\beta$ belong to $\Ux_D$  and $\Lambda_{D,\phi}$ vanish on these elements. Therefore $sX\in \Ux_D$,
$$\Lambda_{D,\phi}(sXs^{-1})=\Lambda_{D,\phi}(sX)=\sum_{\gamma\in D\cap D'}\phi(\gamma)\phi'(\gamma) = c(\phi,\phi')$$
and $\dot{\xi}(sus^{-1}) =  \eps^{c(\phi,\phi')}$. We have $\chi(u) = q^{s(D)-r(D,D')}\eps^{c(\phi,\phi')}$. $\Box$

\section{Supercharacter theory for the group  $G^a$}

The group   $G^a$ is a semidirect product of the subgroup $H$  and the unipotent group  $U^a$.
Recall that the subgroup  $U^a=1+\Ux$ is an algebra subgroup,  and the the group $G^a$ is a group of invertible elements in the associative  algebra   $\gx^a$.

A supercharacter theory for the group  $G^a$ can be constructed following the papers   \cite{P1,P2}.
Consider the group   $\tGa=H\ltimes (U^a\times U^a)$ and its representation in  $\Ux$:
\begin{equation}\label{rhoJ}
\rho(\tau)(X) = tAXB^{-1}t^{-1},
\end{equation}
where  $\tau=(t,A,B)$,~ $X\in \Ux$,~ $t\in H$,~ $A, B\in U^a$.

Define the action of  $t\in H$ on $\phi:D\to\Fq^*$ by the formula  $t.\phi(i,j)=t_i\phi(i,j)t_j^{-1}$.
Theorem  \ref{supclass} implies the classification of   $\tGa$-orbits in $\Ux$: each  $\tGa$-orbit contains some element  $X_{D,\phi}$; two elements   $X_{D,\phi}$ and $X_{D',\phi'}$ belong to a common $\tGa$-orbit if and only if  $D=D'$, and $\phi$ and  $\phi'$ are conjugated with respect to the group $H$ action.

Consider the action of group  $\tGa$ on  $G^a$ by the formula
\begin{equation}\label{rG}
r_\tau(g)=1+ tA(g-1)B^{-1}t^{-1},
\end{equation}
where
$\tau=(t,A,B)\in\tG$.  Here $g-1\in \gx^a$.
The group  $H$ is abelian; it follows the semisimple parts of the elements  $g$ and  $r_\tau(g)$ coincide.

For any rook placement  $D\in\Delta$, we define the subgroup $H_D$, which consists of  all $h=\diag(h_1,\ldots,h_n)$, where $h_i=1$ for  $i\in \row(D)\cup \col(D)$.
\\
{\bf Example}. $n=4$, ~$D=\{\gamma=(2,3)\}$. Then
{\small
	$$ H_D=\left\{\left(\begin{array}{cccc}
	*&0&0&0\\
	0&0&0&0\\
	0&0&0&0\\
	0&0&0&*\\
	\end{array}\right)\right\}.$$}

Consider the set  $\Bc$ that consists of triples  $\bx=(D, \phi, h)$, where $D$ is a rook placement in  $\Delta$, ~$\phi: D\to \Fq^*$, and $h\in H_D$. Attach to  $\bx\in \Bc_0$ the  element
$g_\bx=hu_{D,\phi}$, where $u_{D,\phi}= 1+X_{D,\phi}$.
Define the action of  $t\in H$ on  $\Bc$  by the formula
$t.\bx=(D,t.\phi,h)$.\\
\Theorem \Num\label{Ga-superclass}. 1. Each  $\tGa$-orbit contains  an element  of the form $g_\bx$. 2. Two elements $g_\bx$ and $g_{\bx'}$, where $\bx=(D, \phi, h)$ and  $\bx'=(D', \phi', h')$, belong to a common  $\tGa$-orbit if and only if $\bx$ and  $\bx'$ are conjugated with respect to  $H$.\\
\Proof. It follows from  \cite[Theorem 3.1]{P1}.

We denote by  $K_\bx$ the orbit of  element  $g_\bx$ with respect to the  group  $\tGa$ action.

Consider the set  $\Ac$ that consists of triples  $\ax=(D, \phi, \theta)$, where $D,\phi$ as above, and $\theta$ is a linear character (one dimensional representation) of the subgroup  $H_D$.
Consider the linear character of the subgroup  $G_D=H_DU_D$ defined by the formula
\begin{equation}\label{xixi}
\xi_\ax(g)=\theta(h)\eps^{\la(X)},
\end{equation}
where $g=h+x$, ~ $h\in H_D$,~ $X\in \Ux_D$.
Let $\chi_\ax=\Ind(\xi_\ax,G_D,G^a)$.

Analogously to  the case $\Bx$, one can define the action of $H$ on $\Ac$.
The character $\chi_\ax$ does not depend on action of $H$ on $\phi$ \cite[Proposition 4.2]{P1}.\\
\Theorem \Num \label{supertriang}~\cite[Theorem 4.5]{P1}. \emph{The systems of characters  $\{\chi_\ax\}$ and subsets
$\{K(\bx)\}$, where  $\ax$ (respectively,  $\bx$) runs through the set of representatives of  $H$-orbits in  $\Ac$ (respectively,  $\Bc$), give rise to a supercharacter theory for  $G^a$.}

To calculate supercharacter values on superclasses  we introduce  the following notations:
\\
$$\delta(D,h)=\left\{\begin{array}{ll}
1,~~\mbox{if}~~ h\in H_D,\\
0,~~\mbox{if}~~ h\notin H_D.\\
\end{array}\right.$$
$$\delta(D,h,D') = \delta(D,h)\delta(D,D'),$$
$$R(D,h)=\{(i,j)\in \Sc(D): ~ h_j\ne 1\},\qquad  r(D,h) = |R(D,H)|,$$
$$ m(D,h,D')=s(D)-r(D,D')-r(D,h).$$
For each $t\in H$, denote  $$c_t(\phi,\phi')=\sum_{\gamma\in D\cap D'}\phi(\gamma)\phi'(\gamma)\gamma(t).$$
\Theorem\Num\label{Gvalue}. \emph{The value of supercharacter  $\chi_\ax$, where $\ax=(D,\phi,\theta)$, on the superclass $K_\bx$, where $\bx=(D',\phi',h)$, equals to }
$$\chi_\ax(K_\bx)=\delta(D,h,D')\theta(h)\frac{q^{m(D,h,D')}}{|H_D|}\sum_{t\in H} \eps^{c_t(\phi,\phi')}.$$
\Proof. As a supercharacter is constant on a superclass,~ we have $\chi_\ax(K_\bx) = \chi_\ax(g) $, where $g=hu$,~ $h=(h_1,\ldots,h_n)\in H_{D'}$,~ $u=1+X$ and $X=X_{D',\phi'}$. We get
\begin{equation}\label{chistandard}
\chi_\ax(g) = \frac{1}{|H_D|}\sum \dot{\xi}_\ax(pgp^{-1}),
\end{equation}
where  the sum is taken over all  $p=ts$, ~ $t\in H$, ~~$s\in S(D)$.
Observe that  $pgp^{-1}=h\bmod U^a$. If $h\notin H_D$, then $\chi_\ax(g)=0$.
Similarly to the proof of Theorem  \ref{Uvalue}, if
$D'\cap (\Sc(D)\cup\Sc^\star(D))\ne\varnothing, $
then   $\chi_\ax(g)=0$.\\
\emph{Item 1.} Let further $h\in H_D$  and $D'\cap (\Sc(D)\cup\Sc^\star(D)) = \varnothing$ (i.e. $\delta(D,h,D')\ne 0$). Recall that  $h\in H_{D'}$.  Then  $h\in H_D\cap H_{D'}$ and  $h_i=1$ if $i$ is a row or column number of a root from $D\cup D'$.

Since $s\in S(D)$, the element  $s$ has the form (\ref{SEF}). In this item, we prove that if there exists a root $\al\in R(D,D')\cup R(D,h)$ with  $a_\al\ne 0$, then $pgp^{-1}\notin G_D$ and hence  $\dot{\xi}_\ax(pgp^{-1})= 0$. We prove for  $p=s$ (case $p=ts$ is treated similarly).

Present  $\Sc(D)$  as a union of two subsets
$\Sc(D)=R(D,h)\cup R'(D,h)$, where $R'(D,h)=\{(i,j)\in \Sc(D): ~ h_j= 1\}$.
Observe that  $h_i=1$ for any  $(i,j)\in \Sc(D)$.
The subset   $R(D,D')$ (defined in subsection  \ref{subsecsupU}) is contained in $R'(D,h)$; indeed, if  $\al=(i,j)\in R(D,D')$, then there exists  $\beta=(j,k)\in D'$ and, therefore,   $h_j=1$.

Denote  by  $R_\pm(D,h)$, ~  $R'_\pm(D,h)$  the intersections of relative subsets with $\Delta_+$ and $\Delta_-$.
Consider two subspaces
$$\Yc_1=   \left\{\sum_{\al\in R_+(D,h)} a_\al E_\al + \sum_{\al\in R_-(D,h)}a_\al F_\al\right\},$$
$$\Yc_2=   \left\{\sum_{\al\in R'_+(D,h)} a_\al E_\al + \sum_{\al\in R_-'(D,h)}a_\al F_\al\right\}.$$
We have  $\Yc=\Yc_1\oplus\Yc_2$ and therefore  $\Ux=\Yc\oplus \Ux_D$ and $S(D)=1+\Yc$.

Observe that $\Yc_1\Yc=0$. Really, if   $\al_1=(i,j)\in R(D,h)$,  then $h_j\ne 1 $. Since $h_j=1$ for any root  $\al=(j,k)\in \Sc(D)$,  we get   $E_{\al_1} E_\al=0$  (respectively,  $E_{\al_1} F_\al=0 $ and $F_{\al_1} E_\al=0$).

We make the following remark: for any elements  $Y_1,Z_1\in\Yc_1$ and $Y_2\in \Yc_2$ the equality   $(1+Z_1+Y_2)(1+Y_1+Y_2)^{-1}= 1+Z_1-Y_1$ is true. Indeed,   $(1+Z_1-Y_1)(1+Y_1+Y_2) = 1+Z_1+Y_2$.

The element  $s\in S(D)$ of the form  (\ref{SEF})  can be presented  $s=1+Y_1+Y_2$, there $Y_i\in \Yc_i$.
Denote $s^h=h^{-1}sh$.
We obtain  $$sgs^{-1}= s(h(1+X))s^{-1} = h s^h(1+X)s^{-1}= h(s^hs^{-1}+s^hXs^{-1}).$$

Applying the above remark  we get
$s^hs^{-1} = 1+X_1$, where $X_1\in\Yc_1$  is of the form
$$X_1= \sum_{\al\in R_+(D,h)} a_\al(h_{\col(\al)}-1) E_\al + \sum_{\al\in R_-(D,h)}a_\al (h_{\col(\al)}-1) F_\al. $$
Here $h_{\col(\al)}-1\ne 0$ for each  $\al\in R(D,h)$.   If there exists  $\al\in R(D,h)$ with  $a_\al\ne 0$, then  $X_1$ is a nonzero element from $\Yc_1$.

As  $\Yc_1\Yc=\{0\}$, we get  $X_1s=X_1(1+Y_1+Y_2) = X_1$. Then $X_1s^{-1}=X_1$.
We have $$sgs^{-1} = h( 1+X_1+s^hXs^{-1}) = h(1+(X_1+s^hX)s^{-1}).$$
Let us show that $s^hX$  belongs to  $\Yc_2\oplus\Ux_D$.
Indeed, the elements $s^h$ and   $s$ have the form  (\ref{SEF}).
As in formula (\ref{sX}), the element  $s^hX$ is  a sum of the element  $X$ (under conditions of Item 1 it belongs to  $\Ux_D$) and a linear combination of products  $E_\al E_\beta$,~ $E_\al F_\beta$ and $F_\al E_\beta$, where $\al\in S(D) $ and $\beta\in D'$. Each of these products either equals to  $E_{\al+\beta}$ (respectively, $F_{\al+\beta}$), or it is zero.

If $\al\in R(D,D')$, then the product of matrix units is nonzero for the  only  $\beta\in D'$ and  $\al+\beta\in S(D)$,~ $\col(\al+\beta)=\col(\beta)$. Therefore  $h_{\col(\al+\beta)}=h_{\col(\beta)}=1$. Then  $E_{\al+\beta}$ (respectively,  $F_{\al+\beta}$) belongs to  $\Yc_2$.

If   $\al\notin R(D,D')$, then the products  $E_\al E_\beta$,~ $E_\al F_\beta$ and $F_\al E_\beta$ belong to  $\Ux_D$,  and $\Lambda_{D,\phi}$ vanish on  these products.
We get  $s^hX\in\Yc_2\oplus\Ux_D$. If   $\al\in R(D,D')$ for some  $a_\al\ne 0$, then $s^hX\notin \Ux_D$.

To sum up, we conclude that if there exists  $\al\in R(D,D')\cup R(D,h)$ with  $a_\al\ne 0$, then
$X_1+s^hX\notin \Ux_D$. Since  $\Ux_D$ is a right ideal, we have $(X_1+s^hX)s^{-1}\notin \Ux_D$.
Hence  $sgs^{-1}\notin G_D$ and   $\dot{\xi}_\ax(sgs^{-1})= 0$.  This proves statement of Item 1. \\
\emph{Item  2}.  We finish calculation of  $\chi_\ax(g)$.  It follows from Item 1 that in the formula  (\ref{chistandard}) we may consider  $s$ has the form (\ref{SEF}) and     $a_\al=0$ for each   $\al\in R(D,h)\cup R(D,D')$.  Under this condition   $sgs^{-1}\in G_D$ and $\xi(sgs^{-1})=\theta(h) \eps^{c(\phi,\phi')} $ (see the proof of theorem  \ref{Uvalue}). Respectively, for  $p=ts$, we obtain   $\xi(pgp^{-1})=\theta(h) \eps^{\Lambda_{D,\phi}(tXt^{-1})} = \theta(h) \eps^{c_t(\phi,\phi')} $.
Therefore $$\chi_\ax(g) = \frac{1}{|H_D|}\theta(h) q^{s(D)-r(D,D')-r(D,h)}\sum_{t\in H} \eps^{c_t(\phi,\phi')}. \quad \Box$$

\textit{Authors information:
\\
Aleksandr N. Panov -- apanov@list.ru}

\end{document}